\documentclass[12pt]{article}
\usepackage{amsmath,amsfonts,latexsym}
\numberwithin{equation}{section}

\newcommand{\C}{{\mathbb C}}

\newcommand{\Z}{{\mathbb Z}}
\newcommand{\N}{{\mathbb N}}

\newcommand{\G}{{\cal G}}

\begin{document}
\begin{center}
{\Large{\bf COMPLETE REDUCIBILITY OF INTEGRABLE MODULES FOR THE AFFINE
LIE (SUPER) ALGEBRAS}} \\ [1.7cm]
{\bf S. Eswara Rao} \\
{\bf School of Mathematics} \\
{\bf Tata Institute of Fundamental Research} \\
{\bf Homi Bhabha Road} \\
{\bf Mumbai - 400 005}
 {\bf India} \\ [5mm]
{\bf email: senapati@math.tifr.res.in}
\end{center}

\begin{abstract}
We prove complete reducibility for an integrable module for an affine
Lie algebra where the canonical central element acts non-trivially.
We further prove that integrable modules does not exists for most of
the super affine Lie algebras where the center acts non-trivially.
\footnote{2000 Mathematics subject classification primary
17B65,
17B68 }
\end{abstract}

\pagebreak
\subsection*{Introduction}

Let $\G$ be simple finite dimensional Lie algebra.  Let
$\widehat{\G}$ be the corresponding affine Lie algebra and let $K$ be
the canonical central element.  A module $V$ of $\widehat{\G}$ is
called integrable if the Chevalley generators act locally nilpotently
on $V$.  In [C] the irreducible integrable modules for $\widehat{\G}$
with finite dimensional weight spaces has been classified.  In
particular any irreducible integrable module with finite dimensional
weight spaces where $K$ acts by positive integer is isomorphic to an
highest weight module.  In this work we prove that any integrable
module with finite dimensional weight spaces where $K$ acts by
non-zero scalars is completely reducible (Theorem (1.10)).

The integrable modules where $K$ acts trivially, need not be completely
reducible.  For example consider the $\widehat{\G}$ (without the
derivation) module $\G
\otimes \C [t, t^{-1}]/ (t-1)^2$ where $K$ acts by zero which is not
completely reducible.  (See [E1] for the graded version).

In section 2 we consider affine Lie super algebras and prove that
most often integrable modules with finite dimensional weight spaces
do not exist.  We use stronger definition of the integrability than
that of [KW].  Let $\G$ be simple finite dimensional Lie super
algebra. 
 Let $\widehat{\G}$ be the corresponding affine Lie super algebra.
Assume that it has non-degenerate symmetric invariant billinear form.
Assume that the semisimple part of the even part of $\G$ is \\
at least two components.  Then integrable modules for $\widehat{\G}$ with
finite dimensional weight spaces where center acts by non zero scalar
does not exist.  (Theorem 2.6) Certainly integrable modules with $K$
acting zero exists.  For example loop modules.  Our techniques work only
with the notion of stronger integrability.  We do not know
whether such a result hold with the weaker integrability of
[KW].

In Theorem (2.9), we prove that an integrable irreducible module for
$\widehat{\G}$ with finite dimensional weight spaces where center $K$
acts by positive integer is necessarily a highest weight module,
assuming the semisimple part of the finite even part is only one
component.  In this case we note that (Remark (2.11)) the module is
completely reducible for the even part.  That class includes the affine
Lie super algebras associated with basic Lie super algebras of types
$A(0,n), B (0,n)$ and $C(n)$.

\section*{Section 1}

\paragraph{(1.1)}  We will fix some notations.  All our algebras are
over complex numbers $\C$.  Let $\stackrel{\circ}{\G}$ be simple finite
dimensional Lie algebra.  Let $\stackrel{\circ}{h}$ be a Cartan
subalgebra.  Let $\stackrel{\circ}{Q}$ and $\stackrel{\circ}{\Lambda}$
be root and weight lattice of $\stackrel{\circ}{\G}$.  Let
$\stackrel{\circ}{\Lambda}^+$ be dominant integral weights of
$\stackrel{\circ}{\G}$.   Let $\alpha_1, \cdots, \alpha_n$ be simple
roots and let $\beta$ be highest root of
$\stackrel{\circ}{\G}; \alpha_1^{\vee}, \cdots \alpha_n^{\vee}$ be the
corresponding simple roots. We choose non-degenerate
billinear form on $\stackrel{\circ}{h}^*$ such that $(\beta, \beta)=2$. 

Let $\widehat{\G} = {\G}  \otimes \C [t, t^{-1}]
\oplus \C K \oplus \C d$ be the corresponding untwisted  affine Lie
algbera.
Let $\widehat{h} = \stackrel{\circ}{h}  \oplus \C K \oplus \C d$ be
the Cartan subalgebra of $\widehat{\G}$.  Let $Q$ and $\Lambda$ be the
root and weight lattice of $\widehat{\G}$.  Let $\delta$  be the null
root.  Let $\Lambda_0$ be an element of $\widehat{h}^*$ such that
$\Lambda_0
(\stackrel{\circ}{h})=0, \Lambda_0 (K)=1$ and $\Lambda_0 (d) =0$.  An
element
$\lambda$ in $\stackrel{\circ}{h}^*$ can be treated as an element of
$\widehat{h}^* $ by extending as $\lambda (K) =0$ and $\lambda (d)
=0$. Let $\overline{\lambda}$ be the
restriction  to $\stackrel{\circ}{h}$. Given $\lambda \epsilon
\widehat{h}^*, \lambda$ can be uniquely written as

\paragraph*{(1.2)}  $\lambda = \overline{\lambda}+ \lambda (d)
\delta+\lambda (K) \Lambda_0$.

\paragraph*{(1.3) Definition} An element $\lambda$ in
$\stackrel{\circ}{\Lambda}^+$ is called minimal if for every $\mu
\epsilon \stackrel{\circ}{\Lambda}^+$ such that $\mu \frac{<}{0}
\lambda$ implies $\mu=\lambda$.  Here $\mu \frac{<}{\circ}
\lambda$ means 
$\lambda - \mu = \displaystyle{\sum_{i=1}^{n}} n_i \alpha_i,  n_i
\epsilon \N$. 

\paragraph*{(1.4) Lemma} [H]~ Let $\lambda$ be minimal in
$\stackrel{\circ}{\Lambda}^+$.  Then  $\lambda (\beta^{\vee})=0$ or 1. 

\paragraph*{Proof.}  See Exercise 13 of Chapter III of [H].

\paragraph*{(1.5)}  Let $\lambda \epsilon \Lambda^+$ and
let $V (\lambda)$ be the irreducible integrable highest weight module
for $\widehat{\G}$.  Let $P (\lambda)$ be the set of weights of
$V(\lambda)$.  Define $\mu \leq \lambda$ if $\lambda - \mu
=\displaystyle{\sum_{i=0}^{n}} n_i \alpha_i, n_i \epsilon \N
(\alpha_0$ is the additional simple root of $\widehat{\G})$.  Let
$\overline{P} (\lambda) = \{\overline{\mu} \mid \mu \epsilon P
(\lambda) \}$.  Clearly $\overline{P} (\lambda)$ determines a unique
coset in $\stackrel{\circ}{\Lambda} / \stackrel{\circ}{Q}$.  Let
$\overline{\mu_0}$ be the minimal element in
$\stackrel{\circ}{\Lambda}^+$ in the above coset.  Let $s$ be a
complex number such that $\lambda (d)-s$ is a non-negative integer.

\paragraph*{(1.6) Lemma}  Let $\mu_0= \overline{\mu}_0+ s \delta+
\lambda (K) \Lambda_0$.  Then $\mu_0 \epsilon P (\lambda)$.

\paragraph*{Proof} First note that by minimality of
$\overline{\mu}_0$ we have $\overline{\mu}_0 \frac{<}{0}
\overline{\lambda}$.  Then clearly $\mu_0 \leq \lambda$.

\paragraph*{Claim}  $\mu_0 \epsilon \Lambda^+$.  Consider $\alpha_i^{\vee}$
for  $1 \leq  i \leq n$.  Then $\mu_0 (\alpha_i^{\vee}) = \overline{\mu}_0
(\alpha_i) \epsilon \N$ as  $\overline{\mu}_0 \epsilon
\stackrel{\circ}{\Lambda}^+$.  Now $\alpha_{0}^{\vee} = K - \beta^{\vee}$
and 
$$
\begin{array}{lll}
\mu_0 (\alpha_{0}^{\vee}) &=& \overline{\mu}_0 (K - \beta^{\vee})=
\lambda (K) 
- \overline{\mu}_0 (\beta^{\vee}) 
\end{array}
$$
Now $\lambda (K)$ is positive integer and hence $\lambda (K) \geq 1$.
 We know from Lemma (1.4) that $\overline{\mu}_0 (\beta^{\vee}) = 0 $ or
$1$.  That means 
$$\mu_0 (\alpha_{d+1}^{\vee}) \epsilon \N.$$
This prove the claim.  Thus $\mu_0$ is dominant integral and $\leq
\lambda$.  By Proposition 12.5 (a) of [K1] it follows that $\mu_0
\epsilon P (\lambda)$.

We need the following from [E2].

\paragraph*{(1.7) Lemma}  (Lemma (2.6) of [E2]).  Let $V$ be
integrable module for $\widehat{\G}$ with finite dimensional weight
spaces.  Let $P(V)$ be the set of weights of $V$.  Let $\lambda
\epsilon P(V)$.  Then
\begin{enumerate}
\item[(1)] There exists $\eta_0 \frac{>}{\circ} 0, \eta_0 \epsilon
\stackrel{\circ}{Q}$ such that $\lambda+ \eta_0+\eta \notin
P(V)$ for all $0 \neq \eta \frac{>}{\circ} 0, \eta \epsilon
\stackrel{\circ}{\Lambda}$. 
\item[(2)] There exists $\eta_0^1 \frac{<}{0} 0, \eta^1_0 \epsilon
\stackrel{\circ}{\Lambda}$ such that $\lambda+ \eta_0^1 +\eta \notin
P(V)$ for all $0 \neq \eta \frac{<}{0} 0, \eta \epsilon
\stackrel{\circ}{\Lambda}$
\end{enumerate}

\paragraph*{Proof} (1) follows from the proof of Lemma (2.6) of [E2].
The proof of  (2) is similar.

\paragraph*{(1.8)  Proposition}  Let $V$ be integrable
$\widehat{\G}-$module with finite dimensional weight spaces.  Assume
the canonical central element $K$ acts by positive integers.  Let
$\lambda \in P (V)$.  Then there exists $\eta \geq 0$ such that
$\lambda+ \eta \epsilon \Lambda^+$ and the irreducible integrable
highest weight module $V(\lambda+\eta) \subseteq V$.

\paragraph*{Proof}  By previous lemma there exists $\eta_0
\frac{>}{\circ} 0$ 
such that $\lambda + \eta_0+ \eta \notin P(V)$ for  $0 \neq \eta
\frac{>}{\circ} 0$.  Now by arguments similar to the proof of
Theorem 2.4 (i) of [C] 
will produce an highest weight module with highest weight
$\lambda+\eta_0+\eta_1$ for some $\eta_1 \geq 0$.  Note that
$(\lambda+\eta_0+ \eta_1) (d) \geq \lambda (d)$.   

In the above proof we need our Lemma (1.7) as the proof of Lemma 2.6
(ii) of [C] is 
incomplete.  We now recall the following variation of a standard
result from [K1].

\paragraph*{(1.9) Proposition}  Let $V$ be integrable module for
$\widehat{\G}$ with finite dimensional weight spaces.  Let $K$ act by
positive integer. Suppose for every $v $ in $V$, there exist $N>0$ such
that $U (\widehat{\G})_{\alpha+n \delta} v=0$ for all $n>N$ and for
$\alpha
\epsilon \stackrel{\circ}{\triangle} U \{0\}$.  Then $V$ is completely
reducible. 

We need to recall some standard notations from [K1] and prove two lemmas.

The Cartan subalgebra $h$ carries a non-degenerate billinear form $ ( \
\mid \ )$.
Let $\nu:h \to h^*$ be an isomorphism such that  $\nu (h) (h_1) = (h
\mid h_1)$.  Let $<,>$ be  the induced billinear form on $h^*$.  Recall
the Casimir operator from Section (2.1), $\rho$ in $h^*$ from (2.5) from
[K1].  Also recall the notion of primitive weights from (9.3) of [K1].
Note that in an integrable module the primitive weights are dominant
integral.

\paragraph*{Lemma A}  Let $V$ be as above.  Suppose $\lambda, \mu$ are
primitive weights such that $\lambda - \mu = \beta \in Q^+ - \{0\}$.  Then
$2 < \lambda+ \rho, \nu^{-1} (\beta) > \neq (\beta, \beta)$.

\paragraph*{Proof}  Follows from the proof of theorem (10.7) of [K1]. See
10.7.3 and the next equation in [K1].

\paragraph*{Lemma B} Let $V$ be as above.  Let $v$ be a weight vector of
weight $\lambda$ such that $(\Omega_0 - a I_{V})^k v=0$ for some $k \in
\Z_+$ and $a \in \C$.  Presumably $v^1 \in U_{- \beta} (
\stackrel{\wedge}{\G})v, \beta \in Q$.  Then
$$(\Omega_0 - (a+2 < \lambda+ \rho, \nu^{-1} (\beta) > - (\beta,
\beta) )I_{V}) v^1 =0.$$

\paragraph*{Proof} Follows from (2.6.1) and (3.4.1) of [K1].  Also see
(9.10.2) of [K1].  Note that $V$ is restricted in the sense of [K1].

\paragraph*{Proof of the Proposition} 

Let $\stackrel{\wedge}{\G}= n^- \oplus h \oplus n^+$ be the standard
triangular
decomposition.  Let $V^0 = \{ v \in V \mid n^+v =0\}$.  Clearly $V^0$ is
$h-$ invariant and hence decomposes under $h$.  Let $V^1 = U
(\stackrel{\wedge}{\G}) V^0$.  It is standard fact that in an integrable
module, each highest weight generate an irreducible integrable module.
Thus $V^1$ is completely reducible.  We will now prove that$V=V^1$.

Clearly the Casimir operator $\Omega_0$ acts on $V$ and leaves each finite
dimensional weight space invariant.  Thus $\Omega_0$ is locally finite on
$V$.  Suppose $V \neq V^1$.  Then there exists $v$ in $V -V^1$ of weight 
$\lambda$ such that
$n^+ v \subseteq V^1$ and $(\Omega_0 - a I_V)^k v=0$ for some $k \in \Z_+$
and
$a \in \C$.  Since, clearly $\Omega_0 v \in V^1$, we  have $a=0$ and hence
$\Omega_0^k v=0$.

>From the hypothesis it follows that $U(n^+) v$ is finite dimesnional.  So
it contains vector $u_{\beta} v$ such that $u_{\beta} \in U
(\stackrel{\wedge}{\G})_{\beta}$ and $n^+ u_{\beta} v=0, \beta \in Q^+
-\{0 \}$.  Let $\mu = \lambda+\beta$ and note that $\lambda, \mu$ are
primitive weights.  Thus by Lemma A.
$$2< \mu+ \rho, \nu^- (\beta)> \neq (\beta, \beta). \leqno{(*)}$$

Now by Lemma B it follows that $2< \lambda+\rho, - \nu^-  (\beta)> = (
\beta, \beta)$ as $\Omega_0 (u_{\beta} v)=0$.  This is a contradiction to
$*$.  Thus $V=V^1$ and $V$ is completely reducible.

\paragraph*{(1.10) Theorem} Let $V$ be integrable module
with finite dimensional weight spaces for $\widehat{\G}$.  Suppose all
eigenvalues of $K$ are non-zero.  Then $V$ is completely reducible as
$\widehat{\G}$-module.

\paragraph*{Proof}  First decompose $V$ with $K$ action.  As $K$
commutes with $\widehat{\G}$, each eigenspace is
$\widehat{\G}$-module.  Thus we can assume that $K$ acts by single
scalar.  It is well known that the central element $K$  acts by
integer (see for example [E2]).  Without loss of generality we can assume
that $K$
acts by positive integer.  We now decompose the module
$$V= \oplus_{\lambda  \epsilon \Lambda/Q} \ W_{\lambda}$$
where $\mu_1, \mu_2$ weight occurs in $W_{\lambda}$ then $\mu_1 -
\mu_2 \epsilon Q$.  Clearly each $W_\lambda$ is a
$\widehat{\G}$-module. Thus we can assume that the weights $P(V)$ of
$V$ lie in single coset of $\Lambda$.

\paragraph*{Claim}  Let $\lambda \epsilon P (V)$.  Then there exists
$\eta \geq 0$ such that $\lambda+ \alpha \notin  P(V)$ for all
positive roots $\alpha$ such that $\alpha > \eta$.  

\paragraph*{Proof of the claim}  Suppose there exists infinitely many
positive roots $\alpha$ such that $\lambda + \alpha \epsilon P (V)$.
First by Proposition (1.8) there exists $\eta \geq 0$ such that
$\lambda+\eta \epsilon \Lambda^+, \lambda+ \eta \epsilon P (V),
(\lambda+ \eta) (d) \geq \lambda (d)$ and the irreducible integrable
highest module $V(\lambda+\eta) \subseteq V$.

Let $\overline{P} (V) = \{\overline{\lambda} \mid \lambda \epsilon
{P} (V) 
\}$.  Clearly $\overline{P} (V)$ defines a unique coset in
$\stackrel{\circ}{\Lambda}$.  Let $\overline{\mu_0} $ be the minimal
weight for this coset.  Let $\mu_0= \overline{\mu}_0+ \lambda (d)
\delta+\lambda (K) w \leq \lambda+\eta$.  By lemma (1.6) we have
$\mu_0 \epsilon P (\lambda+\eta) \subseteq P (V)$. 

First note that the number positive roots $\alpha_1$ such that $\alpha_1
\not> \eta$ is finite.

Now choose positive root $\alpha_1 > \eta$ such that $\lambda+
\alpha_1 \epsilon P (V)$.  (This is due to our supposition). Now by
above arguments there exists $\eta_1 \geq 0$ such that $\lambda+
\alpha_1+ \eta_1 \epsilon P (V), \lambda+\alpha_1+\eta_1 \epsilon
\Lambda^+$ and $V(\lambda+\alpha_1+\eta_1) \subseteq V$.  Further
$\mu_0 \leq \lambda+ \alpha_1+\eta_1$ and $\mu_0 \epsilon P (\lambda+
\alpha_1+ \eta_1) \subseteq P (V)$.  Note that $\lambda+ \alpha_1+
\eta_1 > \lambda+ \eta ($ note the strict inequality).  Thus
$V(\lambda+ \alpha_1+ \eta_1) \neq V(\lambda+ \eta)$.  Both modules
have common weight $\mu_0$.  Thus we have proved that dim $V_{\mu_0}
\geq 2$.  By repeating $n$ times the above argument we get dim
$V_{\mu_0} \geq n$.  But dim $V_{\mu_0}$ is finite and thus this
process has to stop.  This proves our claim.
It follows from the claim and that the module $V$ satisfies the conditions
of Proposition 1.9 and hence it is completely reducible.

\paragraph*{(1.11) ~ Remark}  Theorem (1.10) imply that an integrable
module with finite dimensional weight spaces in which $K$ acts by positive
integer belongs to the category ${\cal O}$.

 \section*{Section 2}

Let $\G= \G_0 \oplus \G_1$ be simple finite dimensional Lie super
algbera $\G_0$ (respectively, $\G_1$) being its even (respectively, odd)
part.
We assume that $\G_0$ is reductive.  We further assume that $\G$ carries a
non-degenerate invariant  ``symmetric'' billinear form.  Such Lie super
algebras are
called basic.  We give the list of basis Lie super algberas from
Proposition (1.1) of [K].

$$
\begin{array}{|l|l|l}  \hline
\G& \G_0 \\ [5mm] \hline
A(m,n) & A_{m^+} A_n+ \C,  \\
C (n) & C_n+\C \\
B (m,n) & B_m +C_n \\
D(m,n) & D_m +C_n,  \\
D(2, 1:a) & D_2+A_1 \\
F(4) & B_3 +A_1 \\  
G(3) & G_2+A_1 \\ \hline
\end{array}
$$

In this section we study the integrable representations of the untwisted
affine Lie super algebras of basic Lie super algebras.

Let $\G$ be a basic Lie-super algebra.  Then the restriction to the even
part need not be positive definite.  In fact we choose the form in such a
way that the restriction to the first component of the even part is
positive definite and the restriction to the second component is negative
definite.  (see section 6 of [KW]).  We  normalize the form in such a way
that $(\alpha, \alpha) =2$ where $\alpha$ is the highest root of the first
component of the even part of $\G_0$ and $(\beta, \beta)=-2$ where
$\beta$ is the highest root of the second component. Let $h$ be the Cartan
subalgebra of $\G$  which is contained in the even part.

\paragraph*{(2.1)} Define affine super algebra $\widehat{\G}$. 

$$\widehat{\G}= \G \otimes \C [t, t^{-1}] \oplus \C K \oplus \C d.$$
The Lie bracket is given by the following. Write $x (n) = x \otimes
t^n$. 
$$
\begin{array}{lll}
[x(n), y (m)]&=& [x,y] (m+n)+ n (x,y) \delta_{m+n,0} K \\
\left[ d,x(n) \right] &=& nx(n) ~~ x,y \epsilon \G, m,n \epsilon \Z,
~~ K {\rm is \ central} \\
{\rm Let} \ \widehat{h} &=& h \oplus \C K \oplus \C d
\end{array}
$$
\paragraph*{(2.2) ~~  Definition}  A module $V$ of $\widehat{\G}$ is
called integrable if
\begin{enumerate}
\item[(1)] $V= \displaystyle{\oplus_{\lambda \epsilon \widehat{h}^*}}
V_{\lambda}, V_{\lambda} = \{v \epsilon V \mid h v = \lambda
(h) v, \forall h \epsilon \widehat{h} \}$  
\item[(2)] $V$ is integrable as a $\widehat{\G}_0$ module. 
\item[(3)] For any $v \epsilon V, U (\G) v $ is finite dimensional
\end{enumerate}
Here $U (\G)$ is the universal enveloping algebra of $\G$. 

\paragraph*{(2.3)~~ Remark}

In [KW] integrable modules are studied with weaker condition.  In [KW]
integrability  means the module is integrable only with the affinization
of one simple part of $\G_0$.  Then they have classified irreducible
highest
weight module which are integrable in the above sense.  See Theorem 6.1 \&
6.2 of [KW].

The purpose of this section is to classify irreducible integrable
modules for $\widehat{\G}$ where center $K$ acts non trivially. 

Let $\G_{01}$ and $\G_{02}$ be the first and second and the simple
component of $\G_0$ as above.  (In case of $D (2,n)$ and $D(2,1; \alpha)$
the first component is not simple.  Then we take one of the simple
component). Let $h_1$ and $h_2$ be the respective Cartan sub
algebras. Let $\triangle_1$ and $\triangle_2$ be the corresponding root
system. The following is very standard.  Does not matter whether the
form is positive definite or negative definite.

\paragraph*{(2.4)} For any root $\alpha \epsilon
\stackrel{\circ}{\triangle}_i$, let
$\alpha^{\vee}$ be the co-root.  Let $x_{\alpha}$ be the corresponding
root vector.  Choose $x_{- \alpha}$ in the negative root space such
that $(x_{\alpha}, x_{- \alpha})= \frac{2}{(\alpha, \alpha)}$.  Then
$x_{\alpha}, x_{- \alpha}, \alpha^{\vee}$ is an $s \ell_2$ triple.  Let
$\gamma= \alpha+ n \delta, \alpha \epsilon \stackrel{\circ}{\triangle}_i$.
Let
$\gamma^{\vee} = \alpha^{\vee} + \frac{2n}{(\alpha, \alpha)} K$ be the co-root.
 Then it is easy to check that $x_{\alpha} (n), x_{- \alpha} (-n),
\gamma^{\vee}$ is an $s \ell_2$ triple. 

\paragraph*{2.5 ~~ Lemma}  Let $V$ be an integrable
$\widehat{\G}$-module.  Let $\lambda$ be a weight of $V$.  Let
$\gamma = \alpha+n \delta, \alpha \epsilon \triangle_i$ such that
$\lambda (\gamma^{\vee}) > 0$.  Then $\lambda - \gamma$ is a weight of
$V$. 

\paragraph*{Proof}  Follows from standard $s \ell_2$ theory. 

\paragraph*{(2.6) Theorem}  Notation as above.  Assume the semi
simple part of $\G_0$ has at least two components.  Let $V$ be
integrable module with finite dimensional weight spaces. Let the
central element $K$ act by non-zero scalar.  Then $V$ is necessarily
trivial module.

Without loss of
generality we can assume that $K$ acts by positive integer.  We can
establish the following by the arguments similar to the proof of
theorem (1.10). 

\paragraph*{(2.7)}  For any $\lambda \epsilon P (V)$ there exists
$N>0$ such that 
$$\lambda + \alpha+ n \delta \notin P (V) \ {\rm for \ all} \ n \geq N
\ {\rm and \ for \ all} \ \alpha \epsilon \stackrel{\circ}{\triangle}_1
\cup \{0\}.$$ 

\paragraph*{(2.8)}  There is one problem.  The module $V$ need not have
finite dimensional weight spaces for $\stackrel{\wedge}{\G}_{01} $ as $h_1
\oplus \C K \oplus \C d$ could be much smaller than the Cartan $h=h_1
\oplus h_2 \oplus \C K \oplus \C d$.  To overcome this problem, first
observe that $\stackrel{\wedge}{\G}_{01}$ commutes with $h_2$.  Now
decompose
the module $V$ with respect to $h_2$ and $h_2$ weight space is a
$\stackrel{\wedge}{\G}_{01}$-module with finite dimensional weight spaces.
Now apply arguments similar to the proof of (1.10) to conclude (2.7).

\paragraph*{Claim}  There exists a weight vector $v$ of weight
$\lambda$ such that \\
$x_{\alpha} (n) v=0 \ {\rm for} \ n <0 \ {\rm and \ for \ all} \
\alpha \epsilon \triangle_2 \cup (0)$. 

First we complete the proof assuming the claim.  From the claim we
have $h (n)v=0$ for $n <0$ and $h \epsilon h_2$.  From the standard
Hisenberg highest weight module theory it follows that $h (n) v \neq
0$ for all $n>0$ and for all $h$ in $h_2$.  Thus it follows that
$\lambda+ m \delta$ is a weight for all $m>0$ contradicting (2.7).
Thus the module $V$ has to be trivial.

\paragraph*{Proof of the claim:}  From Lemma 1.7 (2) it follows that there
exists $\lambda \in P (V)$ such that $\lambda - \alpha \notin P (V)$ for
all $\alpha \in \triangle^+_2$.  Let $\triangle^{- a r}_{2}$ be the
negative real roots of $\stackrel{\wedge}{\G}_{02}$.  Define $\triangle
(\lambda) = \{\gamma \in \triangle_2^{-ar}  \mid \lambda (\gamma^{\vee})
 \leq 0 \}$. Then $\triangle (\lambda)$ is finite set.  Indeed, let 
$\gamma = \alpha - n \delta, \alpha \in \stackrel{\circ}{\triangle}_2 
n > 0$ be an element of $\triangle^{-ar}_{2}$.  Then 
$\lambda (\gamma^{\vee}) = \lambda (\alpha^{\vee}) - n
\lambda (K) /_{(\alpha, \alpha)}>0$ for $n$ sufficiently large (recall
$(\alpha, \alpha) <0$ for all $\alpha \in \stackrel{\circ}{\triangle}_2)$.
Fix a positve integer $r$ such that $\alpha - s \delta \in
\triangle^{-a}_{2} - \triangle (\lambda)$ for $s \geq r$.

\paragraph*{Subclaim 1} $\lambda - s \delta \notin P (V)$ for $s \geq
r$.  Suppose $\lambda - s  \delta \epsilon P (V)$ for some $s
\geq r$ we have $\lambda ( (\alpha - s \delta)^{\vee}) > 0$ then by
lemma (2.5), $\lambda - s \delta - (\alpha - s \delta) = \lambda -
\alpha \epsilon P (V)$ which is a contradiction to the choice of
$\lambda$. 

Fix a positive integer $p$ such that $\lambda - s \delta \notin P
(V)$ for $s > p$ and $\lambda - p \delta \epsilon P (V)$.

\paragraph*{Subclaim 2}  $\lambda - \alpha - (m+p) \delta \notin P
(V)$ for $m >0$ and $\alpha \epsilon \stackrel{\circ}{\triangle}^+_2$.
Suppose the
claim is false.  Consider $(\lambda - \alpha - (m + p) \delta)
(\alpha^{\vee}) < 0$ since $\lambda (\alpha^{\vee}) <0$ and $\alpha
(\alpha^{\vee}) =2$.  Then by Lemma (2.5) we have 
$\lambda - \alpha - (m+ p) \delta + \alpha = \lambda - (m+p) \delta
\epsilon P (V)$ contradiction the choice of $p$. 

\paragraph*{subclaim 3} $ \lambda+\alpha - 
(m +p + 1) \delta \notin P (V)
$ for  $m> r$ and $\alpha \epsilon \stackrel{\circ}{\triangle}_2^+$.
Suppose the
claim is false.  Consider

$(\lambda+ \alpha - (m+p+1)
\delta)(\alpha - m \delta)^{\vee} >0$ as $\alpha - m \delta \notin
\triangle (\lambda)$.  Thus by lemma 2.5 we have
$$\lambda+ \alpha - (m+1+p) \delta - \alpha +m \delta$$
$$=\lambda - (1+p) \delta \ \epsilon P (V)$$
contradicts the choice of $p$.  

Thus we have proved

$$\stackrel{\wedge}{\G}_{02, -r \delta} V_{\lambda - p \delta} =0, {r>0} \
{\rm and} \ \
\stackrel{\wedge}{\G}_{02, \alpha - s \delta} V_{\lambda - p \delta} =0 \
{\rm for \ all}$$
but finitely many negative roots.  Since $V$ is integrable $\overline{W}=
U(\widehat{\G}_{02}^{-}) V_{\lambda - p \delta}$ is finite dimensional.
Let
$\mu$ be the lowest weight of $\overline{W}$.  This weight satisfies all
the requirements of the claim.

\paragraph*{(2.9) ~~ Theorem}  Let $\widehat{\G}$ be the affine super
algebra defined earlier.  Assume that the semisimple part of the finite
even part has only one component.  Further assume that the
non-degenerate form restricted to this simple Lie algebra $\G_0$ is
positive definite.  Let $V$ be irreducible integrable module with
finite dimensional weight spaces.  Assume the central element $K$
acts as positive integer.  Then $V$ is an highest weight module. 

\paragraph*{Proof}  From the proof of Theorem (2.6) we have (2.7).
Let $\beta_1, \cdots, \beta_k$ be odd roots of $\G$.  Let $v$ be a
wieght vector of $V$ of weight $\lambda$. 

\paragraph*{Claim}  The following vectors span is a finite dimensional
space $W$
$$\{ x_{\beta_{i_1}} (m_1) \cdots x_{\beta_{i_k}} (m_{k}) v, i_j \leq
i_k, m \geq 0 \}$$
where $x_{\beta_i}$ is a root vector for the odd root space
$\G_{\pm \beta_i}$. The affine roots that are occuring in the product
are all distinct. In the above we take negative roots first and postive
roots next.  The indices are decreasing order.
It is sufficient to prove that the vector space $T$ spanned by the
following vector is finite dimensional.
$$\{x_{\beta} (m_1) \cdots x_{\beta} (m_k) w, m_i \geq 0, m_i \neq m_j 
\}$$
where $w$ is any weight vector of  $V$.  
This is because there are only finitely  many odd roots in $\G$.

First note that if $k \beta$ is a root for $k>0$ then   $k=1$ or $2$.
Consider $x_{\beta} (m) w = [h (m), x_{\beta}] w = h (m) x_{\beta} (w)
\pm x_{\beta}  h (m) w$.  By (2.7) both vectors are zero for large $m$.
Let $m_0$ be such that $x_{k \beta} (m) w=0$ for $m > m_0$ and $k=1,2$.
Then it is easy to see that $T$ is spanned by
$$\{x_{\beta} (m_1) \cdots x_{\beta} (m_k)w; 0 \leq m_i \leq  m_0,  m_i
\neq m_j,  i \neq j\}$$
which is clearly finite dimensional.

Let $H$ be the center of the reductive Lie algebra $\G_0$.  Consider
$$S=U (\widehat{\G}_0^-) U (h)U (\widehat{\G}_0^+) U
(\displaystyle{\oplus_{n>0}} H \otimes t^n) W. \leqno(2.10)$$ 
 $${\rm Then} \ V=U (\displaystyle{\oplus_{n<0}} H \otimes t^n)U
(\widehat{\G}_1^-)S$$
by PBW basis theorem. 

By (2.7) we conclude that $U (\displaystyle{ \oplus_{n>0}} H \otimes t^n)
W=W_1$ is finite dimensional.  Clearly $S$ is $\widehat{\G}_0$-module and
by Theorem (1.10)  $S$ is  completely reducible.  In fact it is direct sum
of highest weight modules.  Since $W_1$ is finite
dimensional it intersects only finitely many of them.  Say $V(\lambda_1) 
\cdots V(\lambda_k)$.  Thus $S=\oplus V (\lambda_i)$ a finite sum.
Thus $S$ has a maximal weight. (Here the ordering is the following $\mu_1
\leq \mu_2$ means $\mu_2 - \mu_1 = \sum n_i \alpha_i,  n_i \in \N,
\alpha_i'$s are small roots of $\widehat{\G}$. $\widehat{\G}$ is a
generalized Kac-Moody Lie super algebra and it does admit simple roots.
See [KW]).  The maximal weight is in fact maximal for $V$ as the rest of
the space brings the weights down.

The maximal weight is in fact highest weight.  As $V$ is
irreducible, it is irreducible highest weight module.

\paragraph*{(2.11) Remark}  In the process we also established that
an irreducible integrable highest weight module for $\widehat{\G}$ is
completely reducible for $\widehat{\G_0 \oplus H}$. 

\paragraph*{Proof}  Let $V$be irreducible highest weight module for
$\widehat{\G}$.  Let 
$$\Omega (V) = \{v \epsilon V \mid h (k) v =0 \ {\rm for \ all} \ h
\epsilon H, \ k>0\}$$
Let $M(k)$ be the irreducible highest weight module for $H \otimes \C
[t, t^{-1}] \oplus \C K$  where $K$ acts by $k$.  Then by Theorem
(1.7.3) of [FLM] we have $V= \Omega (V) \otimes M (k)$.  Now $\Omega
(V)$ is an integrable module and hence by Theorem (1.10) decomposes
into irreducible modules for $\widehat{\G}_0$.  Thus the Remark
follows. 

\pagebreak

\begin{center}
{\bf REFERENCES}
\end{center}
\vskip 1cm

\begin{enumerate}
\item[{[C]}] Chari, V. {\it Integrable representations of Affine
Lie-Algebras.}  Invent Math. 85, 317-335 (1986).
\item[{[E1]}] Eswara Rao, S. {\it Classification of Loop Modules with
finite dimensional weight spaces}. Math. Ann. 305, 651-663 (1996).
\item[{[E2]}]  Eswara Rao, S. {\it Classification of irreducible
integrable modules for Toroidal Lie algebras with finite dimensional
weight spaces } preprint 2001.
\item[{[FLM]}] Frenkel, I., Lepowsky, J. and  Mueurman, A. {\it Vertex
operator algebras and the Monster}, Academic Press (1989).
\item[{[H]}] Humpreys, J.E., {\it Introduction to Lie-algebras and 
representation theory}, Springer Berlin, Hidelberg, New York (1972).
\item[{[K]}] Kac, V.G., {\it Representations of classical Lie super
algebras} Lecture note in Math. 676 (1978), 597-626. 
\item[{[K1]}] Kac, V. {\it Infinite dimensional Lie-algebras}, Cambridge University Press, 3rd edition (1990).
\item[{[KW]}] Kac, V.G. and Minoru Wakimoto {\it Integrable highest
weight modules over affine superalgebras and appell's function},
Communications of Mathematical Physics, 215, 631-682 (2001).
\end{enumerate}
\end{document}